\documentclass[12pt]{article}
\usepackage{epsfig,amssymb,amsmath}

\usepackage{graphicx}

\newcommand{\sta}{{\rm sta}}

\newcommand{ \PP}{{P}}

\newcommand{\vsig}{\varsigma}

\newcommand{\barvsig}{\bar{\varsigma}}

\newcommand{\barx}{\bar{x}}

\newcommand{\real}{{\mathbb R}} 

\newcommand{\half}{\frac{1}{2}}

\newcommand{\bvsig}{\mbox{\boldmath$\varsigma$}}

\newcommand{\veps}{\xi}

\newcommand{\bveps}{\mbox{\boldmath$\xi$}}

\newcommand{\barbvsig}{\bar{\mbox{\boldmath$\varsigma$}}}

\newcommand{\eba}{\begin{array}}
\newcommand{\eea}{\end{array}}
\newcommand{\ebe}{\begin{eqnarray}}
\newcommand{\eee}{\end{eqnarray}}
\newcommand{\eb}{\begin{equation}}
\newcommand{\ee}{\end{equation}}

\newcommand{\calP}{{\cal{P}}}

\newcommand{\bx}{{\bf x}}

\newcommand{\calS}{{\cal S}}

\newcommand{\calE}{{\cal E}}

\newcommand{\calX}{{\cal X}}

\newcommand{\barbx}{\bar{\bf x}}

\newcommand{\alp}{{\alpha}}

\newcommand{ \Lam}{{\Lambda}}

\newcommand{ \xx}{{\bf x}}

\newtheorem{thm}{Theorem}

\newcommand\VV{V}

\renewcommand\Pi{{{P}}}

\renewcommand\xx{{{x}}}

\renewcommand\barbx{{\bar{x}}}

\renewcommand\eb{\begin{equation}}
\renewcommand\ee{\end{equation}}

           \newcommand\PPd{{P^d}}

\renewcommand\xx{{{\bf x}}}
\renewcommand\barbx{{\bar{\bf x}}}

\renewcommand\barbx{{\bar{\xx}}}

\newcommand\delshap{\delta^\sharp}
\newcommand\delflet{\delta^\flat}

\begin{document}
\noindent {\bf {\Large Canonical Duality Theory for Solving Minimization Problem of Rosenbrock Function}\\ [0.5cm]
David Yang Gao $\cdot$ Jiapu Zhang}\\ [0.5cm]
{\em Graduate School of Information Technology and Mathematical Sciences,\\
University of Ballarat, Victoria, Australia\\[0.2cm]
}

\noindent {\bf Abstract} This paper presents a {\em canonical duality theory } for solving
nonconvex minimization problem of Rosenbrock function. Extensive numerical results show
that this benchmark test problem can be solved precisely and efficiently  to obtain 
global 
optimal solutions.\\


\section{Introduction}
\begin{figure}[h!]
\centerline{
\includegraphics[scale=0.47]{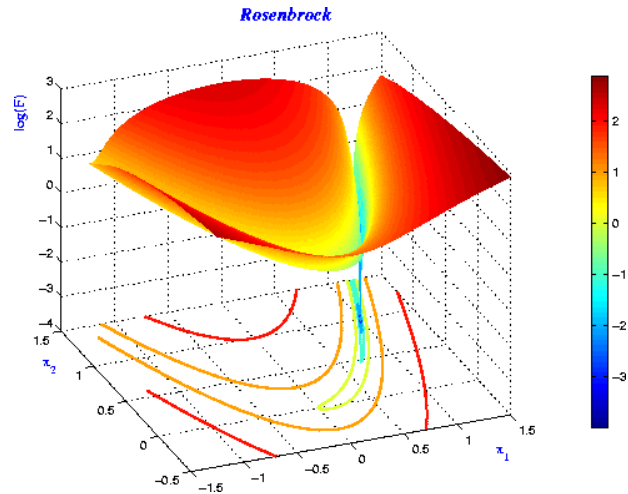}
} \caption{2-dimensional Rosenbrock function (www2.imm.dtu.dk/courses/02610/)} \label{2d_rosenbrock_function}
\end{figure}

Nonconvex minimization problem
of Rosenbrock function,  introduced in \cite{rosen60}, is a benchmark test problem
in global optimization
that has been used extensively to test   performance of optimization algorithms and numerical
approaches. The global minimizer of this function
is located in a long, deep, narrow, parabolic/banana shaped flat valley (Figure \ref{2d_rosenbrock_function}).

Although  to find this  valley is trivial for most cases,
 to accurately locate the global optimal solution is
   very difficult by almost all gradient-type methods and some derivative-free methods.
Due to the nonconvexity, it can be easily tested that
   if the initial point is chosen to be  $(3,3, \dots, 3)$,
   direct algorithms  are always trapped into
   a local minimizer for problems with dimensions $n = 5 \thicksim 7 $ as well as
   $n  \ge  4000$;   if the initial point  is chosen
    at $(100,100,\dots,100)$,   iterations will be stopped at a local min
    with the objective function value
     $> 47.23824896$ even for a two-dimensional problem.
   This paper will show that by the canonical duality theory,   this well-known
   benchmark problem can be solved efficiently in an elegant way.\\

The {\em canonical duality theory} was originally developed in nonconvex/nonsmooth mechanics
 \cite{gao-dual00}. It is now  realized that this   potentially powerful theory can be used for
 solving a large class of nonconvex/nonsmooth/discrete problems \cite{gao-optm03, gao2010}.
   In this short research note,  we will first show the
    nonconvex minimization problem of Rosenbrock function can be
    reformulated as a canonical dual problem (with zero duality gap)
    and the critical point of the Rosenbrock function  can be analytically expressed in terms of
    its canonical dual solutions.
  Both global and local extremal solutions can be identified
by the triality theorem.
Extensive numerical examples and discussion are presented in the
  last section.

\section{Primal Problem and Its Canonical Dual}
 The primal problem is
\eb
    (\calP):\;\;\;   \min  \left\{ P (\xx) =
    \sum_{i=1}^{n-1}  \left[ (x_i - 1)^2 + \half \alp (x_{i+1} - x_i^2)^2 \right] \;\;
    | \;\; {\xx \in \calX } \right\}
  , \label{eq:govp}
  \ee
 where  $\xx = \{ x_i \} \in \calX = \real^n$  is a real unknown vector,
  $\alp =2 N$ and $N$  is a
 given real number.
 Clearly, this is a nonconvex minimization problem which could have multiple local minimizers.\\

 In order to use the canonical duality theory for solving this nonconvex problem, we need
 to define a {\em geometrically admissible} canonical measure
 \eb
 \bveps = \{ \veps_i\} = \{ x_i^2 -  x_{i+1} \} \in \calE_a \subset \real^{n-1}.
 \ee
 The canonical function $\VV:
 \calE_a \rightarrow \real$ can be defined by
 \eb
 \VV(\bveps) = \sum_{j=1}^{n-1} \half \alp \veps_j^2 ,
 \ee
 which is a convex function. The canonical dual variable $\bvsig = \bveps^*$ can be defined
 uniquely by
 \eb
 \bvsig = \{ \vsig_j \} = \nabla \VV (\bveps) = \{ \alp \veps_j \} .
 \ee
 Therefore, by the Legendre transformation, the conjugate function $\VV^*:\calS =\real^{n-1}
 \rightarrow \real$ is obtained as
 \eb
 \VV^*(\bvsig) = \sta \{ \bveps^T \bvsig - \VV(\bveps) \; | \; \bveps \in \calE_a \}
 = \sum_{j=1}^{n-1} \half \alp^{-1} \vsig_j^2 .
 \ee
 Replacing $\sum_{i=1}^{n-1}   \half \alp (x_{i+1} - x_i^2)^2 $ by
 the Legendre equality $ \VV(\Lam(\bx)) = \Lam(\bx)^T \bvsig - \VV^*(\bvsig)$,
 the total complementary function $\Xi:\calX \times \calS \rightarrow \real$ is
 given by
 \eb
 \Xi(\bx, \bvsig) = \sum_{i=1}^{n-1}  \left[ (x_i - 1)^2 +  \vsig_i (x_i^2 - x_{i+1}  )
 - \half \alp^{-1} \vsig_i^2   \right]
 . \label{complementary_function}
 \ee\\

Let $\delshap $ and $\delflet$ be shifting operators such that
$\delshap \vsig_i = \vsig_{i+1}$ and $\delflet \vsig_i = \vsig_{i-1}$.
We define $\delflet \vsig_1 = 0$.
Then on the canonical dual feasible space
\eb
\calS_a = \{ \bvsig \in \calS | \;\; \vsig_i + 1 \neq 0 \;\; \forall
i = 1, \dots, n-2, \; \vsig_{n-1} = 0 \}, \label{dual_feasible_space}
\ee
the canonical dual can be obtained by
\eb
\PPd(\bvsig) = \sta \{ \Xi(\bx , \bvsig) | \; \bx \in \calX \}
= n-1 - \sum_{i=1}^{n-1} \left[
\frac{(\delflet \vsig_i + 2 )^2}{4 (\vsig_i + 1)} + \half \alp^{-1} \vsig_i^2 \right]. \label{dual_problem1}
\ee
Based on the {\em complementary-dual principle} proposed in  the canonical duality theory
(see \cite{gao-cc2009}),
we have the following result.

\begin{thm} If $\barbvsig$ is a critical point of $\PPd(\bvsig)$, then
the vector $\barbx = \{\barx_i\} $ defined by
\eb
\barx_i = \frac{\delflet \barvsig_i + 2 }{2 (\barvsig_i + 1)},  \;\; i=1, \cdots, n-1,
\;\;   \barx_{n}=\barx^2_{n-1}
 \label{eq-anasol}
\ee
is a critical point of $\PP(\bx)$ and
\eb
\PP(\barbx) = \Xi(\barbx, \barbvsig) = \PPd(\barbvsig).
\ee
\end{thm}
This theorem presents actually an ``analytic" solution  form to the Rosenbrock function,
i.e. the critical point of the Rosenbrock function must be in the form of (\ref{eq-anasol})
for each dual solution $\barvsig$.
The first version of this analytical solution form was presented in
nonconvex variational problems in
phase transitions and finite deformation mechanics \cite{gao97,gaomecc99,gaona00}. 
The extremality of the analytical solution is governed by the so-called {\em triality theory}.
Let
\begin{eqnarray}
 \calS_a^+ = \{ \bvsig \in \calS_a | \;\; \vsig_i + 1 > 0 \;\; \forall i = 1, \dots , n-1\},
 \end{eqnarray}
 we have the following theorem:
 \begin{thm}
 Suppose that $\barbvsig$ is a critical point of $\PPd(\bvsig)$ and
the vector $\barbx = \{\barx_i\} $ is defined by Theorem 1.

If $\barbvsig \in \calS_a^+ $,
then $\barbvsig$ is a global maximal solution to the canonical dual problem on $\calS^+_a$, i.e.,
\eb
(\calP^d_+): \;\; \max \{ \PPd(\bvsig) \; | \; \bvsig \in \calS_a^+ \},  \label{equ-dual}
\ee
the vector  $\barbx $ is a global minimal to the primal problem, and
\eb
\PP(\barbx) = \min_{\bx \in \calX } \PP(\bx) = \max_{\bvsig \in \calS_a^+ } \PPd(\bvsig)
= \PPd(\barbvsig).
\ee
\end{thm}

\noindent Theorem 2 shows that the canonical dual problem $(\calP^d_+)$
provides a global optimal solution to the nonconvex primal problem.
Since $(\calP^d_+)$ is a concave maximization problem over a convex space
which can be solved easily.
This theorem is actually a special application of Gao and Strang's general result
on global minimizer in in nonconvex analysis
\cite{gao-gs89a}.\\

By introducing
\eb
\calS_a^- = \calS / \calS_a^+=\{ \bvsig \in  \real^{n-1} | \;\; \vsig_i + 1 < 0 \;\; \forall i = 1, \dots , n-1\} ,
\ee
recently the triality theory (see \cite{gao-wu-jogo10}) leads to the following theorem.
\begin{thm}
Suppose that $\barbvsig$ is a critical point of $\PPd(\bvsig)$ and
the vector $\barbx = \{\barx_i\} $ is defined by Theorem 1.

If $\barbvsig \in \calS_a^- $,
then on a neighborhood $\calX_o \times \calS_o \subset \calX \times \calS_a^-$
of $(\barbx, \barbvsig )$,
we have either
\eb
\PP(\barbx) = \min_{\bx \in \calX_o} \PP(\bx) =
 \min_{\bvsig \in \calS_o}   \PPd(\bvsig) = \PPd(\barbvsig)   , \label{equ-dualmin}
\ee
or
\eb
\PP(\barbx) = \max_{\bx \in \calX_o} \PP(\bx) =
 \max_{\bvsig \in \calS_o}   \PPd(\bvsig) = \PPd(\barbvsig)   , \label{equ-dualmax}
\ee
\end{thm}

\noindent The proof of this Theorem can be derived from the recent paper by Gao and Wu \cite{gao-wu-jogo10}.
By the fact that the canonical dual function is a d.c. function
(difference of convex functions) on $ \calS_a^-$,
 the double-min duality (\ref{equ-dualmin}) can be used for finding the
biggest
 local minimizer of the Rosenbrock function $\PP(\bx)$,
 while the double-max duality (\ref{equ-dualmax})
  can be used for finding the biggest local maximizer of
  $\PP(\bx)$.
In physics and material sciences, this pair of biggest local extremal points play important roles
 in phase transitions.\\
 
Because $\vsig_{n-1} = 0$, we may know that $ \calS_a^-$ is an empty set. 
Thus, by Theorem 3 in this paper we cannot find a local maximizer or minimizer on $ \calS_a^-$ or its subset for $\PPd(\barbvsig)$.

\section{Numerical Examples and Discussion}
$(\mathcal{P})$ and $(\mathcal{P}^d_+)$ will be solved by the discrete gradient (DG) method (\cite{Adil_DiscreteGradientMethod}), which is a local search optimization solver for nonconvex and/or nonsmooth optimization problems. In two dimensional space, Rosenbrock function has a long ravine with very steep walls and flat bottom; ``because of the curved flat valley the optimization is zig-zagging slowly with small stepsizes towards the minimum" (en.wikipedia.org/wiki/Gradient\_descent). This means any gradient method may fail to minimize the Rosenbrock function even from 2 dimensions. The DG method is a derivative-free method which can be applied for miminizing/maximizing Rosenbrock function and its dual. Numerical experiments have been carried out in Intel(R) Celeron(R) CPU 900@2.20GHz Windows Vista$^{\text{TM}}$ Home Basic personal notebook computer.\\

We try $N$=100 (when $N$=10 we find the numerical results are similar to $N=100$), with the dimensions $n$=2$\thicksim$10, 20, 30, 40, 50, 60, 70, 80, 90, 100, 200, 300, 400, 500, 600, 700, 800, 900, 1000, 2000, 3000, 4000. We first set $(3,3,\dots,3)$ (called seed1) as the initial solution for ($\mathcal{P}$) (usually the feasible solution space is a box constrained between -2.048 and 2.048 \cite{ai_etal1998, xin_etal2009, yu_and_wang2008}). Numerical results (Table 1) show that to solve
 the primal problem ($\mathcal{P}$),  the DG method can easily and quickly get approximate global minimum solution  to $\barbx = (1, 1, \dots, 1)$ with the  approximate global optimal values at
  $P(\barbx) = 0,$ except for $n$=5$\thicksim$7, 4000, where the DG method can only get a  local minimum solution $\barbx = (-1, 1, \dots, 1)$ with $P(\barbx) = 4$. Then we let
   $\bx_0 = (100,100,\dots, 100)$ (called seed2) be the initial solution for ($\mathcal{P}$),
    searched in the intervals $-500\leq x_i\leq 500, i=1,2,\dots,n$. We find that the
    DG method was trapped into local optimal solutions but not getting any global minimum at all, even from a 2 dimensional problem (see Table 2), its objective function value is 47.23824896. However, from Table 2 we can  see that by the same DG method,
 the global maximum of the dual problem can be obtained very elegantly.\\

For $(\mathcal{P}^d_+)$, the corresponding dimensions are 1$\thicksim$9, 19, 29, 39, 49, 59, 69, 79, 89, 99, 199, 299, 399, 499, 599, 699, 799, 899, 999, 1999, 2999, 3999.
The initial solution is set as $\bvsig_0= (-2/3, -2/3, \dots,-2/3,0)$ (called seed1),
the constraints $\bvsig_i +1 >0, i=1,2,\dots,n-1$ were penalized into the objective function;
by $\Xi (\bx, \bvsig  )'_{x_n}=0$ of formula (6),  we can set the values of the last variable $\varsigma_{n-1}$ always being 0 ($>-1$). With these numerical computation settings, the DG method can easily and quickly solve all these test problems  to accurately get a global maximizer
$\barbvsig = (0,0,...,0)$  with the optimal  value $\PP^d(\barbvsig) = 0$ (Table 1).
By the fact that the canonical dual problem
$(\mathcal{P}^d_+)$ is a concave maximization over a convex open space,
the  DG method was not trapped into any local optimal solution.
But, for the nonconvex primal problem ($\mathcal{P}$) in dimensions $n$=5$\thicksim$7 and 4000,
 the DG method was trapped into  local minimizer  $\barbx = (-1, 1, \dots, 1)$.
 If we set the initial solution as $\bvsig_0 = (100, 100, \dots,100,0)$ (called seed2)
 and repeat the calculations, our numerical results (Table 2) show again that the
 canonical dual problem can be
   easily and quickly solved by the DG method
    to accurately get the global maximizer
    $\barbvsig = (0,0,...,0) $ with the optimal solution $\PP^d(\barbvsig) = 0$
        for dimensions $n= 1  \thicksim 9, 19, 29, 39, 49, 59, 69, 79, 89, 99, 199, 299, 399, 499, 599, 699, 799, 899, 999, 1999$.\\

The comparisons between $(\mathcal{P})$ and $(\mathcal{P}^d_+)$ in view of total number of iterations and total number of objective function evaluations (i.e. function calls) are listed in Tables \ref{table1}-\ref{table2}. Compared with $(\mathcal{P}^d_+)$, the approximate global and local optimal solutions and their optimal objective function values of $(\mathcal{P})$ are not accurate, and even cannot be obtained  if  the  initial iteration is set to be
 $\bx_0 = (100,100, \dots, 100)$. In Table \ref{table1}, we can see that the total number of iterations and function calls for  $(\mathcal{P})$ are always greater than those for
  $(\mathcal{P}^d_+)$. This means that $(\mathcal{P}^d_+)$ costs less computer calculations than $(\mathcal{P})$, though $(\mathcal{P}^d_+)$ still can get accurate global optimal solutions and the global optimal objective function value.
The initial solutions $\bx_0 = (100,100, \dots, 100)$ and $\bvsig_0 =
(100,100, \dots, 100, 0)$ respectively for $(\mathcal{P})$ and $(\mathcal{P}^d_+)$ are not practical for real  numerical tests so that the total number of iterations and function calls of $(\mathcal{P})$ are sometimes less than those of $(\mathcal{P}^d_+)$. Regarding the CPU times for  solving
$(\mathcal{P}^d_+)$ with
$n=4000$,  the largest CPU time for seed1 is 206.3581 seconds (i.e. 3.4393 minutes).\\

\begin{table}[h!]
\begin{center}
\caption{Results of numerical experiments for $(\mathcal{P})$ and $(\mathcal{P}^d_+)$: $N=100$, seed1}
\label{table1} {\small
\begin{tabular}{|c|c|c|c|c|c|c|} \hline
Dimension $n$ &\multicolumn{2}{|c|}{Iterations}    &\multicolumn{2}{|c|}{Function calls}    &\multicolumn{2}{|c|}{Objective function value}\\ \cline{2-7}
    &$(\mathcal{P})$ &$(\mathcal{P}^d_+)$  &$(\mathcal{P})$ &$(\mathcal{P}^d_+)$      &$(\mathcal{P})$ &$(\mathcal{P}^d_+)$\\ \hline
2   &120             &24                 &2843            &28                     &0.00001073&0.00000000\\ \hline
3   &422             &26                 &8996            &137                    &0.00401438&0.00000000\\ \hline
4   &3737            &35                 &48352           &202                    &0.00615273&0.00000000\\ \hline
5*  &335             &34                 &10179           &399                    &3.96077434&0.00000000\\ \hline
6*  &2375            &44                 &43770           &868                    &4.00635895&0.00000000\\ \hline
7*  &1223            &53                 &28009           &1625                   &4.09419146&0.00000000\\ \hline
8   &2160            &55                 &46792           &2100                   &0.01246714&0.00000000\\ \hline
9   &2692            &51                 &61017           &2526                   &0.01397307&0.00000000\\ \hline
10  &4444            &63                 &91470           &3979                   &0.01055630&0.00000000\\ \hline
20  &3042            &55                 &140924          &10084                  &0.00940077&0.00000000\\ \hline
30  &2321            &58                 &133980          &20515                  &0.01075478&0.00000000\\ \hline
40  &1659            &60                 &173795          &26818                  &0.01227866&0.00000000\\ \hline
50  &2032            &57                 &219233          &36459                  &0.01264147&0.00000000\\ \hline
60  &1966            &61                 &260701          &50495                  &0.01048188&0.00000000\\ \hline
70  &1876            &56                 &272919          &52545                  &0.01531147&0.00000000\\ \hline
80  &1405            &61                 &195156          &59684                  &0.01594730&0.00000000\\ \hline
90  &2142            &61                 &371963          &71320                  &0.01055831&0.00000000\\ \hline
100 &2676            &60                 &510722          &70208                  &0.01125514&0.00000000\\ \hline
200 &1395            &61                 &653604          &188589                 &0.01115318&0.00000000\\ \hline
300 &1368            &60                 &882760          &235163                 &0.01574873&0.00000000\\ \hline
400 &2085            &66                 &1869675         &301805                 &0.00928066&0.00000000\\ \hline
500 &1155            &59                 &1394240         &358938                 &0.01168440&0.00000000\\ \hline
600 &1226            &63                 &1808285         &451817                 &0.00918730&0.00000000\\ \hline
700 &1557            &60                 &2134359         &559378                 &0.01257100&0.00000000\\ \hline
800 &1398            &61                 &2098062         &522726                 &0.01442714&0.00000000\\ \hline
900 &716             &65                 &1904187         &763449                 &0.01074534&0.00000000\\ \hline
1000&1825            &61                 &3598608         &681509                 &0.00897202&0.00000000\\ \hline
2000&257             &62                 &2087277         &1455472                &0.00937219&0.00000000\\ \hline
3000&3221            &60                 &20642543        &2714296                &0.01250373&0.00000000\\ \hline
4000*&679            &60                 &7581502         &3659292                &4.11193171&0.00000000\\ \hline
\end{tabular}
}
\end{center}
\end{table}

\begin{table}[h!]
\begin{center}
\caption{Results of numerical experiments for $(\mathcal{P})$ and $(\mathcal{P}^d_+)$: $N=100$, seed2}
\label{table2} {\small
\begin{tabular}{|c|c|c|c|c|c|c|} \hline
Dimension $n$ &\multicolumn{2}{|c|}{Iterations}    &\multicolumn{2}{|c|}{Function calls}    &\multicolumn{2}{|c|}{Objective function value}\\ \cline{2-7}
    &$(\mathcal{P})$ &$(\mathcal{P}^d_+)$  &$(\mathcal{P})$ &$(\mathcal{P}^d_+)$      &$(\mathcal{P})$ &$(\mathcal{P}^d_+)$\\ \hline
2   &10013           &24                   &227521          &28                       &47.23824896     &0.00000000\\ \hline
3   &144             &32                   &4869            &235                      &96.49814330     &0.00000000\\ \hline
4   &144             &81                   &5279            &938                      &82.46230602     &0.00000000\\ \hline
5   &148             &137                  &5682            &1768                     &94.19254867     &0.00000000\\ \hline
6   &154             &166*                 &6238            &2590                     &88.84382963     &0.00000000\\ \hline
7   &159             &179*                 &7097            &3288                     &237.63078399    &0.00000000\\ \hline
8   &165             &202*                 &7502            &4300                     &238.41126013    &0.00000000\\ \hline
9   &153             &206*                 &7137            &5083                     &84.54205412     &0.00000000\\ \hline
10  &162             &216*                 &7491            &5920                     &83.23094398     &0.00000000\\ \hline
20  &225             &285*                 &19111           &17458                    &83.94779152     &0.00000000\\ \hline
30  &216             &301*                 &20939           &28543*                   &156.95838274    &0.00000000\\ \hline
40  &163             &291*                 &19775           &40444*                   &83.30960344     &0.00000000\\ \hline
50  &158             &298*                 &33269           &51888*                   &85.93091895     &0.00000000\\ \hline
60  &158             &312*                 &34094           &61767*                   &89.07412094     &0.00000000\\ \hline
70  &162             &284*                 &35436           &69865*                   &92.45725362     &0.00000000\\ \hline
80  &209             &297*                 &35607           &89127*                   &157.69955825    &0.00000000\\ \hline
90  &227             &294*                 &60398           &98748*                   &82.44035053     &0.00000000\\ \hline
100 &202             &290*                 &57792           &102796*                  &81.94595276     &0.00000000\\ \hline
200 &1826            &262                  &436413          &189293                   &83.77165551     &0.00000000\\ \hline
300 &195             &259*                 &169238          &261320*                  &152.95671738    &0.00000000\\ \hline
400 &195             &278*                 &212104          &375816*                  &82.49253919     &0.00000000\\ \hline
500 &190             &297*                 &331637          &522695*                  &82.40170647     &0.00000000\\ \hline
600 &292             &303*                 &431092          &559068*                  &150.15456693    &0.00000000\\ \hline
700 &189             &275*                 &383735          &758631*                  &89.14575473     &0.00000000\\ \hline
800 &198             &270*                 &429674          &701053*                  &84.50538257     &0.00000000\\ \hline
900 &198             &280*                 &416150          &867398*                  &85.32757049     &0.00000000\\ \hline
1000&193             &283*                 &445326          &930761*                  &89.48369379     &0.00000000\\ \hline
2000&232             &310*                 &1123240         &2030104*                 &84.26810981     &0.00000000\\ \hline
\end{tabular}
}
\end{center}
\end{table}

{\bf Example 1.} Let $n=4$ (four dimensions). The global minimizer
 is known to be $\barbx = (1,1,1,1)$ and $\PP(\barbx) = 0$.\\
\noindent {\bf Solution:} By using the DG method for both primal problem  $(\mathcal{P})$
and its canonical dual $(\calP^d_+)$, we have the
numerical solutions
\[
\barbx =(1.0166873133, 1.0337174892, 1.0687306765, 1.1425101552), \;\;\; P(\barbx) =0.00615273,
 \]
\[
\barbvsig =(0.0000000119, 0.0000000000, 0.0000000000), \;\;
P^d_+ (\barbvsig)=0.00000000.
\]
 This shows that  the canonical dual problem provides more accurate solution.   \\

{\bf Example 2.} For dimension $n=5$, the Rosenbrock function has exactly two minima,
one is the global optimal solution $(1,1,1,1,1)$ with global optimal minimum value $0$, 
and another minimum is a local minimum near $(-1,1,1,1,1)$ with local optimal minimum value 4.\\
\noindent {\bf Solution:} By the DG method,
the primal solution is
\[
\barbx =(-0.9856129203, 0.9814803343, 0.9682775584, 0.9398661046, 0.8840549028)
\]
with $P(\barbx) =3.96077434. $
Clearly, this is a local minimizer.
While the canonical dual problem produces accurately a global optimal solution
\[
\barbvsig =(0.0000004388, 0.0000006036, 0.0000000000, 0.0000000000), \;\;
P^d_+ (\barbvsig ) =0.
\]

{\bf Example 3.} For $n=6$ (six dimensions), the Rosenbrock
function has exactly two minima, i.e., the  global optimal solution
\[
\barbx_1 = (1,1,1,1,1,1), \;\; \PP(\barbx_1) = 0,
\]
and local minimal solution
\[
\barbx_2 = (-1,1,1,1,1,1), \;\; \PP(\barbx_2) =  4.
\]

\noindent {\bf Solution:} To solve the primal problem directly, the DG method can only provide
local solution
\[
\barbx =(-0.9970726441, 1.0041582933, 1.0133158817, 1.0292928527, 1.0607123926, 1.1258344785)
\]
with $  P(\barbx) =4.00635895 $.
For the canonical dual problem, the DG method produces
\[
\barbvsig =(0.0000001747, -0.0000000559, 0.0000005919, 0.0000000000, 0.0000000000),\]
\[
P^d_+ (\barbvsig) =0.
\]

{\bf Example 4.} Similarly, if  $n=7$, the test problem has  the same
global optimal  solution
\[
\barbx_1 = (1,1,1,1,1,1,1), \;\; P(\barbx_1) = 0
 \]
 and  the local minimal solution
 \[
 \barbx_2 = (-1,1,1,1,1,1,1) , \;\; \PP(\barbx_2 ) =  4.
 \]

\noindent {\bf Solution:} By the DG method, we have
\begin{eqnarray*}
\barbx &=& (-1.0003403494, 1.0106728675, 1.0264433859, 1.0561180077, \\
 & & 1.1168007274, 1.2483026410, 1.5594822181) ,\\
  P(\barbx) &= & 4.09419146 ,
\end{eqnarray*}
\begin{eqnarray*}
\barbvsig  &= & (-0.0000001431, -0.0000011147, -0.0000010643, -0.0000003284,\\
& &  0.0000000000, 0.0000000000),\\
P^d_+ (\barbvsig) & = & 0.
\end{eqnarray*}
This shows again that the DG iterations for solving the primal problem is
trapped to a local min.

{\bf Example 5.} Now we let $n=4000$.
The Rosenbrock function has many minima. The global
 optimal solution is still
 $\barbx_1 = (1,\dots,1)$ with  $P(\barbx) = 0$.
 One of local minima is  nearby the point
 $\barbx_2 = (-1,1,\dots,1)$ with  $P(\barbx_2) =4.$\\

\noindent {\bf Solution:} Again, by the DG method, the primal iteration
is trapped at
\[
\barbx =(-0.9932861006, 0.9966510741,  \dots, 1.3122885708, 1.7233744896), \;\;
P(\barbx) =4.11193171.
\]
The conical dual solution is
\begin{eqnarray*}
\barbvsig &=& (-0.0000000314, -0.0000000040, -0.0000000437,  \dots , \\
& & -0.0000000281, 0.0000000008, -0.0000000214, 0.0000000000, 0.0000000000),
\end{eqnarray*}
which produce precisely the optimal value
$P^d_+ (\barbvsig) =0. $
 Indeed,  as long as  $n\geq 5$,
   the DG method is always trapped into the local minimizer  $\barbx = (-1, 1, \dots , 1)$
   if the initial solution is set to be
 $\bx_0 = (-1.0005, 1.0005, \dots , 1.0005)$.  \\

It is worth to note that both $P(\xx)$  and $\PPd(\bvsig)$
 are the sum of $n-1$ items. This is convenient for MPI (Message Passing Interface) parallel computations. We may broadcast (MPI\_Bcast) the sum to $n-1$ processes, each process calculates one item, and at last all the partials are reduced (MPI\_Reduce) onto one process to get the sum. Thus on Tambo machines of VLSCI (http://www.vlsci.unimelb.edu.au) we should be able to successfully solve (\ref{eq:govp}) and (\ref{equ-dual}) with at least $3.2767\times 10^7$ variables if setting the maximal variables for the DG method to be 4000 (though the DG method and its parallelization version (\cite{gleb_etal}) can solve optimization problems with more than 4000 variables). The successfully tested MPI code is followed:
{\sf
\begin{enumerate}
\item[] broadcast $n-1$
      \begin{enumerate}
      \item[] call MPI\_BCAST ($n-1$,1,MPI\_INTEGER,  0, MPI\_COMM\_WORLD , ierr)
      \end{enumerate}
\item[] check for quit signal
       \begin{enumerate}
      \item[] if ( $n-1$ .le. 0 ) goto 30
      \end{enumerate}
\item[] calculate every partials
       \begin{enumerate}
      \item[] sum  = 0.0d0
      \item[] do 20 i = myid$+1$, $n-1$, numprocs
                \begin{enumerate}
                 \item[] sum $=$ sum $+ (x(i)-1.0)**2 +100.0*(x(i)**2 -x(i+1))**2$
                 \end{enumerate}
      \item[] 20 continue (\underline{for $P (\xx)$})
      \item[] do 20 i = myid$+1$, $n-1$, numprocs
                \begin{enumerate}
                 \item[] if ($i-1$ .eq. 0) then $\vsig (0)$=0
                 \item[] sum $=$ sum $+ (\vsig (i-1)+2.0)/(4*(\vsig (i)+1.0)) +(1.0/400.0)*\vsig (i)**2$
                 \end{enumerate}
      \item[] 20 continue (\underline{for $\PPd(\bvsig)$})
      \item[] f $=$ sum
      \end{enumerate}
\item[] collect all the partial sums
      \begin{enumerate}
      \item[] call MPI\_REDUCE (f,objf,1,MPI\_DOUBLE\_PRECISION, MPI\_SUM, 0,\\
     \& MPI\_COMM\_WORLD,  ierr )
      \end{enumerate}
\item[] 30 node 0 (i.e. myid = 0) prints the sums = objf
\end{enumerate}
}

\section{Conclusion}
This research note demonstrates a powerful application of the
 {\em canonical duality theory } for solving the nonconvex minimization problem of Rosenbrock function. Extensive numerical computations show that by using the same DG method,
 the canonical dual problem can be easily solved to produce
  global solutions.

\section*{{\small Acknowledgments:}}
This research is supported by US Air Force
Office of Scientific Research under the grant AFOSR FA9550-10-1-0487. The authors thank Professor Kok Lay Teo (Curtin University, Australia) for his helpful comments and Associate Professor Adil M. Bagirov (Ballarat University, Australia) for his FORTRAN codes of his discrete gradient method.

 \end{document}